\documentclass[11pt]{article}
\usepackage{fullpage}

\usepackage{graphicx} 

\usepackage[a4paper, total={5.5in, 8.5in}]{geometry}

\usepackage{amsmath,amsthm,amssymb}
\usepackage[T1]{fontenc} 
\usepackage[utf8]{inputenc} 
\usepackage{graphicx}
\usepackage{scrextend}
\usepackage{subcaption}

\usepackage{xcolor}

\usepackage{hyperref} 

\newtheorem{theorem}{Theorem}

\newtheorem{problem}{Problem}
\newtheorem{definition}{Definition}

\newtheorem{prop}[theorem]{Proposition}

\begin{document}

\title{\textsc{Max-Cut} on Token Graphs}

\author{
Jos\'e Angel Alejandro-Soto\thanks{Departamento de Matem\'aticas, CINVESTAV,
Mexico. \texttt{jaalejandro@math.cinvestav.mx},
\texttt{ruyfabila@math.cinvestav.edu.mx}, \texttt{sgomez@math.cinvestav.mx}}
\and
Paulo T. Araujo\thanks{Department of Physics and Astronomy, University of Alabama, Tuscaloosa, AL, 35487.
\texttt{ptaraujo@ua.edu}, \texttt{mashamim@crimson.ua.edu}}
\and
Ruy Fabila-Monroy\footnotemark[1]
\and
Sergio Gerardo G\'omez-Galicia\footnotemark[1]
\and
Mohamed Hibat-Allah\thanks{Department of Applied Mathematics, University of
Waterloo, Ontario N2L 3G1, Canada; and Vector Institute, Toronto, Ontario
M5G 0C6, Canada. \texttt{mhibatallah@uwaterloo.ca}}
\and
Md Moshiur Rahman Raj\thanks{Department of Physics, University of Rajshahi, Rajshahi 6205, Bangladesh. \texttt{moshiur\_rahman@tutanota.com}}
\and
Mahmud Ashraf Shamim\footnotemark[2]
\and
Ana Trujillo-Negrete\thanks{Instituto de Matem\'aticas, UNAM, Mexico.
\texttt{atrujillo@im.unam.mx}}
}


\maketitle

\begin{abstract}
Let $G$ be a graph on $n$ vertices, with a positive weight assigned to each edge.
A \emph{cut} of $G$ is a partition of its vertices into two sets $X$ and $Y$; the \emph{weight} of $(X,Y)$ is the sum of the weights of the edges of $G$ with one endpoint in $X$ and the other endpoint in $Y$. The \textsc{Max-Cut} problem asks for the maximum weight, $m(G)$, of a cut of $G$. The \emph{$k$-token graph} of $G$ is the graph, $F_k(G)$,  whose vertices are the $k$-subsets of vertices of $G$, two of them being adjacent when their symmetric difference is an edge of $G$. The phase-reconstruction problem for Heisenberg antiferromagnets in physics can be rephrased as a \textsc{Max-Cut} problem on token graphs. Motivated by this, we study the \textsc{Max-Cut} problem on token graphs, mainly when all edges are assigned
the same weight. Our main focus is to bound $m(F_k(G))$ in terms of $m(G)$. We prove that $m(F_k(G))\geq\binom{n-2}{k-1}m(G)$, and that this bound is tight up to a multiplicative factor of $1+O(k/n)$; in the case when $k=c \cdot n$ for some constant $0 < c < 1$, we provide an example where this bound is off by a multiplicative constant factor.
We also show bounds for $m(F_k(G))$ for various classes of graphs $G$. Finally, we consider the algorithmic complexity of finding $m(F_k(G))$, showing it to be NP-complete for a constant $k$, and giving evidence that, when the weights take finitely many values, it is not NP-hard for $k\geq n^{\delta}$ for any constant $\delta>0$.
\end{abstract}

\section{Introduction}

Let $G=(V,E)$ be a simple graph on $n$ vertices. A \emph{cut} $(X,Y)$ of $G$ is a partition of its vertices into two sets $X$ and $Y$. Its \emph{cut set} is
\[
E(X,Y):=\{xy\in E:x\in X,y\in Y\}.
\]
Suppose that each edge $e$ of $G$ is assigned a nonnegative weight $w(e)$, and let
$w(G):=\sum_{e\in E(G)}w(e)$. The \emph{weight} of $(X,Y)$ is defined as
\[
w(X,Y):=\sum_{xy \in E(X,Y)} w(x,y).
\]
A \emph{maximum cut} of $G$ is a cut of maximum weight; we denote this weight by
\[
m(G):=\max\{w(X,Y):(X,Y)\text{ is a cut of }G\}.
\]
\textsc{Max-Cut} is the algorithmic problem of computing a maximum cut of a given graph. In the \emph{unweighted} case, in which every edge of $G$ has weight equal to one, we have that $w(X,Y)=|E(X,Y)|$; and we may also refer to the weight of $(X,Y)$ by its size. \textsc{Max-Cut} is an important problem
in combinatorial optimization; it arises naturally in numerous application domains, including VLSI design, circuit layout and statistical physics~\cite{Barahona1988, GraphPartitioningSurvey}.

The decision problem of \textsc{Max-Cut}, that is, the problem of determining whether $G$ has a cut
of weight at least a given value, was shown to be NP-complete by Karp~\cite{KarpList} in its weighted version; the unweighted case was settled shortly after by Garey, Johnson, and Stockmeyer~\cite{GJS76}.
\textsc{Max-Cut} is also hard to approximate.
A uniformly random partition cuts half of the edges in expectation; therefore, in general
\[m(G) \ge \frac{1}{2}w(G).\]
This yields a $1/2$-approximation; the algorithm of Goemans and Williamson~\cite{GoemansWilliamson1995} achieves an approximation ratio of $\approx 0.87856$ using semidefinite programming and randomized rounding. Assuming the Unique Games Conjecture, this ratio is optimal~\cite{KKMO07}. It is also known that  it is NP-hard to approximate \textsc{Max-Cut} within a factor greater than $16/17\approx 0.941$~\cite{Hastad01,TSSW00}, even in the unweighted case~\cite{crescenzi2001weighted}.

Recently, applications of \textsc{Max-Cut} in quantum many-body
physics~\cite{PhysicMathLanguage, Shamim:2026hiv, apte2026} have been found.
In the Heisenberg antiferromagnet, a quantum spin sits at each vertex of $G$
(the physical lattice) and interacts along its edges. Each spin can be ``up'' or ``down''; and
the Hamiltonian of the system conserves the number $k$ of ``up'' spins. A state with exactly $k$ up
spins is described by the set of $k$ vertices at which the up spins sit; the up spins may thus be
regarded as $k$ indistinguishable \emph{tokens} placed on distinct vertices of $G$, and the
Hamiltonian acts by moving one token at a time along an edge of $G$ to an unoccupied vertex.

We can represent this behavior with the following combinatorial construction. Let $1 \le k \le n-1$.
The \emph{$k$-token graph} of $G$ is the graph $F_k(G)$ whose vertices are all the $k$-subsets
of vertices of $G$, with $A$  being adjacent to $B$ if their symmetric difference $A \triangle B$ is an edge of $G$.
We refer to $G$ as the \emph{parent graph} of $F_k(G)$.
Token graphs have been introduced independently several times and under different names: as \emph{$k$-subgraph graph}~\cite{Johns}, as
\emph{double vertex graphs} in the case $k=2$~\cite{Alavi1991}, as \emph{$n$-tuple vertex
graphs}~\cite{Zhu1992},  as \emph{symmetric powers} of graphs~\cite{Rudolph2002, Audenaert2007} and token graphs~\cite{Token}.
Throughout this paper we follow the notation and terminology of~\cite{Token}. The states with $k$ up spins correspond to the vertices of
$F_k(G)$, and the transitions allowed by the Hamiltonian to its edges~\cite{PhysicMathLanguage}.

Since this Hamiltonian is a real symmetric matrix, the ground state of the system can be written with real
coefficients $s_A\psi_A$, where $\psi_A\geq 0$ and $s_A\in\{+1,-1\}$ for each vertex $A$ of $F_k(G)$.
Reconstructing the signs $s_A$ when the amplitudes $\psi_A$ are prescribed---the \emph{phase
reconstruction problem}---amounts to maximizing a weighted cut of $F_k(G)$, in which each edge $AB$
receives the weight $\psi_A\psi_B$ scaled by the coupling strength of the edge of $G$ along which the
token moves~\cite{PhysicMathLanguage}.
When $G$ is bipartite, $F_k(G)$ is bipartite as well, the
maximum cut consists of all edges, and the optimal signs are given by the Marshall sign
rule~\cite{marshall1955}. When $G$ is not bipartite (\emph{frustrated}, in physical terminology), odd
cycles in $F_k(G)$ make the problem a nontrivial \textsc{Max-Cut} instance. The physically most
relevant case is when the numbers of up and down spins coincide, corresponding to $k=n/2$;
in this case, the ground
state of an antiferromagnet is expected to be a total spin singlet~\cite{LiebMattis1962}, which lies
in this sector.

This quantum setting contrasts with the classical connection between
\textsc{Max-Cut} and statistical physics, in which computing a ground state
of the Ising model amounts to computing a maximum cut of the physical
interaction graph itself~\cite{Barahona1988, Galluccio2001}; once the
dynamics are restricted to a fixed number of up spins, the optimization
problem is naturally defined not on the interaction graph, but on its token
graph.
This distinction is significant, as restricting the problem to token graphs
may lead to stronger structural results than are available for general graphs.

In this paper we study the \textsc{Max-Cut} problem on token graphs, both from a combinatorial and an algorithmic perspective. Throughout, every edge $AB$ of $F_k(G)$ inherits the weight of the edge $A\triangle B$ of $G$. When these weights are the coupling strengths, the Hamiltonian restricted to the states with $k$ up spins equals a constant minus half the Laplacian of $F_k(G)$~\cite{Shamim:2026hiv}. In the phase reconstruction problem the weight of an edge $AB$ is the coupling strength of $A\triangle B$ multiplied by $\psi_A\psi_B$; because of this factor, the weighting of $F_k(G)$ is in general not induced by a weighting of $G$, and we comment on this case only in Section~\ref{sec:complexity}. In Section~\ref{sec:parents} we restrict ourselves to the unweighted case; for the nearest-neighbor antiferromagnet this is the case of equal coupling strengths.

Noting that $A\mapsto V(G)\setminus A$ is an isomorphism between
$F_k(G)$ and $F_{n-k}(G)$, we assume throughout that $k\leq n/2$. Finally, for a
positive integer $r$ we write $[r]:=\{1,\dots,r\}$, and $\binom{[r]}{j}$ for the
family of $j$-subsets of $[r]$.

\section{Bounds for $m(F_k(G))$ in terms of $m(G)$}\label{sec:bounds}

Let $(X,Y)$ be a cut of $G$. From  $(X,Y)$ we can obtain a cut of $F_k(G)$ as follows. Let
\[
X':=\{A\in V(F_k(G)): |A\cap X|\text{ is even}\}
\]
and
\[
Y':=\{A\in V(F_k(G)): |A\cap X|\text{ is odd}\}.
\]
We call $(X',Y')$ the \emph{parity lift} of $(X,Y)$.

Let $S$ be a set of $k-1$ vertices of $G$, and let $xy$ be an edge of $G$ such that
$x,y \notin S$.
Note that $xy\in E(X,Y)$ if and only if $S\cup\{x\}$ and $S\cup\{y\}$ lie on different sides of $(X',Y')$. That is,  the edges of $E(X',Y')$
 are  exactly the edges of $F_k(G)$ whose symmetric difference lies
in $E(X,Y)$; as every edge $xy$ of $G$ is the symmetric difference of
$\binom{n-2}{k-1}$ edges of $F_k(G)$, each of weight $w(xy)$, the parity lift has
weight
\[
    w(X',Y')
    =
    \binom{n-2}{k-1}w(X,Y).
\]
Thus,
\begin{equation}\label{eq:lower}
    m(F_k(G))\geq \binom{n-2}{k-1}m(G).
\end{equation}

This bound is tight when $G$ is bipartite, but as we now show it is not tight in
general, even in the unweighted case.

\begin{prop}\label{prop_K5_local}
\[
    20=m(F_2(K_5))> \binom{5-2}{2-1}m(K_5)=18.
\]
\end{prop}
\begin{proof}
 For $i\in [5]$, let $H_i$ be the subgraph of $F_2(K_5)$ induced by
\[
    \{\{i,j\}:j\in[5],\ j\neq i\}.
\]
The edge sets of $H_1,\ldots,H_5$  are a partition of $E(F_2(K_5))$.
Since each $H_i$ is a complete graph on $4$ vertices and every cut of $K_4$ contains
at most $4$ edges we have that every cut of
$F_2(K_5)$ has size at most $5\cdot 4=20$;
this is attained with the cut
\[
   ( \{12, 23, 34, 45, 15\},\{13,14,24,25,35\}).
\]
\end{proof}


We use the previous example to construct arbitrarily large graphs where the lower bound of (\ref{eq:lower}), for a linear value of $k$,
is off by a constant multiplicative factor greater than one.
\begin{theorem}\label{thm:K5q}
Let $G$ be the disjoint union
of $q$ copies, $H_1,\dots,H_q$ of $K_5$; let $n=5q$ and  $0<\alpha\leq 1/2$ be a constant, such that $k=\alpha n$ is an integer.  Then
\[
    m(F_k(G))
    \geq
    \left(1+\frac{\alpha(1-\alpha)}{3}-o(1)\right)
    \binom{n-2}{k-1}m(G).
\]
In particular, for $\alpha=1/2$ and $q$ even we have that
\[
    m(F_{n/2}(G))
    \geq
    \left(\frac{13}{12}-o(1)\right)
    \binom{n-2}{n/2-1}m(G).
\]
\end{theorem}
\begin{proof}
We begin by noting that
\[m(G)=6q.\]
For every $1 \le i \le q$ and
$1\leq r\leq 4$, fix a maximum cut $(X_{i,r},Y_{i,r})$ of $F_r(H_i)$.
For every
$A \in V(F_k(G))$ and every $1 \le i \le q$, let
\[\sigma_i(A):=\begin{cases} +1 & \textrm{ if } A\cap V(H_i) \in X_{i,r}
\text{ for some } 1 \le r \le 4; \\[1ex]
-1 & \text{ if } A\cap V(H_i) \in Y_{i,r}
\text{ for some } 1 \le r \le 4;\\[1ex]
+1 & \text{ if }  A\cap V(H_i)=\emptyset \text{ or }  A\cap V(H_i)=V(H_i).
\end{cases}\]
Let
\[X':=\left \{ A\in V(F_k(G)): \prod_{i=1}^q \sigma_i(A)=+1 \right \}
\text{ and }
Y':=\left \{ A\in V(F_k(G)): \prod_{i=1}^q \sigma_i(A)=-1 \right \}.\]
Note that $(X',Y)$ is a cut of $F_k(G)$.

Since there are no edges between the $H_i$, for every edge $AB \in E(X',Y')$, we have that
\[\sigma_i(A)\neq \sigma_i(B),\]
for some $1 \le i \le q$, and
\[\sigma_j(A)= \sigma_j(B)\]
for all other $j \neq i$.
Thus, for a fixed $H_i$ and a fixed $r$, each edge of $(X_{i,r},Y_{i,r})$
can be extended to an edge of $(X',Y')$ by choosing the remaining $k-r$ tokens
outside $H_i$, in $\binom{5(q-1)}{k-r}$ ways.  This implies that
\[
m(F_k(G))\geq
q\sum_{r=1}^4 m(F_r(K_5))\binom{5(q-1)}{k-r}.
\]
Note that
\[
m(F_1(K_5))=m(F_4(K_5))=m(K_5)=6, \text{ and } m(F_2(K_5))=m(F_3(K_5))=20,
\]
where the last equalities follow from
Proposition~\ref{prop_K5_local} and the fact that
$F_k(G) \simeq F_{n-k}(G)$.

Therefore,
\[
m(F_k(G))\geq
q\left[
6\binom{5q-5}{k-1}
+20\binom{5q-5}{k-2}
+20\binom{5q-5}{k-3}
+6\binom{5q-5}{k-4}
\right].
\]
For every $1 \le r \le 4$, we have
\[
\frac{\binom{5q-5}{k-r}}{\binom{5q-2}{k-1}}
=
\frac{\prod_{j=1}^{r-1}(k-j)\,\prod_{j=1}^{4-r}(5q-k-j)}{(5q-2)(5q-3)(5q-4)}.
\]

Since $k=\alpha 5q$, as $q \to \infty$ we have that $\prod_{j=1}^{r-1}(k-j) \sim (\alpha 5q)^{r-1}$, $\prod_{j=1}^{4-r}(5q-k-j) \sim ((1-\alpha)5q)^{4-r}$, and $(5q-2)(5q-3)(5q-4) \sim (5q)^3$.
Thus, as $q \to \infty$,
\[
    \frac{\binom{5q-5}{k-r}}{\binom{5q-2}{k-1}}
    \longrightarrow
    \alpha^{r-1}(1-\alpha)^{4-r}.
\]
Therefore,
\[
\begin{aligned}
\frac{m(F_k(G))}
{\binom{5q-2}{k-1}m(G)}
&\geq
\frac{
6(1-\alpha)^3
+20\alpha(1-\alpha)^2
+20\alpha^2(1-\alpha)
+6\alpha^3
}{6}
-o(1)\\
&=
1+\frac{\alpha(1-\alpha)}{3}-o(1);
\end{aligned}
\]
the result follows.
\end{proof}

Let $S \in \binom{V(G)}{k-1}$ and $(X',Y')$ be a cut of $F_k(G)$. We show an upper bound on $m(F_k(G))$
by projecting $(X',Y')$ to a cut of $G-S$;  where $G-S$ denotes the subgraph of $G$ induced by $V(G)\setminus S$. Let
\[
X_S:=\{x\in V(G)\setminus S:S\cup\{x\}\in X'\}
\]
and
\[
Y_S:=\{y\in V(G)\setminus S:S\cup\{y\}\in Y'\}.
\]

Let $AB \in E(F_k(G))$ and $\{x,y\}:=A\triangle B$.
Note that $AB\in E(X',Y')$ if and only if $xy\in E(X_{A\cap B},Y_{A\cap B})$.
Therefore,
\begin{equation}\label{eq:cut_grouping}
w(X',Y')=\sum_{S\in\binom{V(G)}{k-1}} w(X_S,Y_S),
\end{equation}
and
\begin{equation}\label{eq:cut_grouping_ub}
m(F_k(G))\leq\sum_{S\in\binom{V(G)}{k-1}} m(G-S).
\end{equation}
We now use inequality (\ref{eq:cut_grouping_ub}) to give an upper bound
on $m(F_k(G))$ in terms of $m(G)$, $n$ and $k$.
\begin{theorem}\label{thm:upper_bound}
\[
m(F_k(G))
\leq
\frac{
2\binom{n}{k-1}\binom{n-2}{k-1}
}{
\binom{n}{k-1}+\binom{n-2}{k-1}
}
m(G).
\]
In particular:
\begin{enumerate}
\item[$a)$] if $k=o(n)$, then
\[
m(F_k(G))
=
\left(1+O\!\left(\frac{k}{n}\right)\right)\binom{n-2}{k-1}m(G);
\]
\item[$b)$] if $k=\alpha n+o(n)$, for some constant $0\leq \alpha\leq 1/2$, then
\[
m(F_k(G))
\leq
\left(
\frac{2}{1+(1-\alpha)^2}
+o(1)
\right)
\binom{n-2}{k-1}m(G);
\]
and
\item[$c)$] if $n$ is even and $k=n/2$, then
\[
m(F_{n/2}(G))
\leq
\left(\frac85+o(1)\right)
\binom{n-2}{n/2-1}m(G).
\]
\end{enumerate}
\end{theorem}
\begin{proof}
Take a maximum cut of $G-S$ and extend it to a cut of $G$ by
placing every vertex of $S$ independently and uniformly at random on one
of the two sides. Every edge not in $G-S$ is cut with probability $1/2$.
Therefore,
\[
m(G)\geq m(G-S)+\frac{w(G)-w(G-S)}{2};
\]
which implies
\[
m(G-S)
\leq
m(G)-\frac{w(G)-w(G-S)}{2}.
\]
By (\ref{eq:cut_grouping_ub}), we have
\[
\begin{aligned}
m(F_k(G))
&\le \sum_{S\in \binom{V(G)}{k-1}} m(G-S)\\
&\leq
\sum_{S\in \binom{V(G)}{k-1}}
\left(
m(G)-\frac{w(G)-w(G-S)}{2}
\right) \\
&=
\binom{n}{k-1}m(G)
-\frac12\binom{n}{k-1}w(G)
+\frac12
\sum_{S\in \binom{V(G)}{k-1}} w(G-S) \\
&=
\binom{n}{k-1}m(G)
-\frac12
\left(
\binom{n}{k-1}-\binom{n-2}{k-1}
\right)
w(G).
\end{aligned}
\]
Where the last equality follows from the fact that every edge $xy \in E(G)$ and every $S \subset V(G) \setminus \{x,y\}$ gives
way to exactly one edge of $F_k(G)$.
Since $m(F_k(G))\le w(F_k(G))=\binom{n-2}{k-1}w(G)$,
\begin{equation}\label{eq:minbound}
m(F_k(G))
\leq
\min\left\{
\binom{n-2}{k-1}w(G),\;
\binom{n}{k-1}m(G)
-\frac12
\left(
\binom{n}{k-1}-\binom{n-2}{k-1}
\right)
w(G)
\right\}.
\end{equation}
 We regard both expressions as functions of $w(G)$.
The first expression is increasing in $w(G)$, while the second is decreasing. Therefore,
the minimum of the two expressions is at
most the value at which they coincide. That is, when
\[
\binom{n-2}{k-1}w(G)
=
\binom{n}{k-1}m(G)
-\frac12
\left(
\binom{n}{k-1}-\binom{n-2}{k-1}
\right)
w(G).
\]
Solving this equation gives
\[
w(G)
=
\frac{2\binom{n}{k-1}}
{\binom{n}{k-1}+\binom{n-2}{k-1}}
m(G).
\]
Substituting this value of $w(G)$ into
$\binom{n-2}{k-1}w(G)$, we obtain
\[
m(F_k(G))
\leq
\frac{
2\binom{n}{k-1}\binom{n-2}{k-1}
}{
\binom{n}{k-1}+\binom{n-2}{k-1}
}
m(G).
\]

We now show $a)$, $b)$ and $c)$. Since
\[
\binom{n-2}{k-1}
=
\frac{(n-k+1)(n-k)}{n(n-1)}\cdot \binom{n}{k-1},
\]
we have that
\[
\frac{
2\binom{n}{k-1}\binom{n-2}{k-1}
}{
\binom{n}{k-1}+\binom{n-2}{k-1}
}
=
\frac{2n(n-1)}
{n(n-1)+(n-k+1)(n-k)}
\binom{n-2}{k-1}.
\]

If $k=o(n)$, then
\[
\frac{2n(n-1)}
{n(n-1)+(n-k+1)(n-k)}
-1
=
\frac{(k-1)(2n-k)}
{n(n-1)+(n-k+1)(n-k)}
\leq
\frac{2kn}{n(n-1)}
=
O\!\left(\frac{k}{n}\right),
\]
and therefore, combining with the lower bound (\ref{eq:lower}), we have that
\[
m(F_k(G))
=
\binom{n-2}{k-1}m(G)
\left(1+O\!\left(\frac{k}{n}\right)\right);
\]
which proves $a)$.

If $k=\alpha n+o(n)$ with $0\leq\alpha\leq 1/2$, then
\[
\frac{(n-k+1)(n-k)}{n(n-1)}=(1-\alpha)^2+o(1).
\]
Therefore,
\[
\frac{2n(n-1)}{n(n-1)+(n-k+1)(n-k)}
=\frac{2}{1+(1-\alpha)^2+o(1)}
=\frac{2}{1+(1-\alpha)^2}+o(1);
\]
which proves $b)$.

If $n$ is even and $k=n/2$, then
\[
\frac{2n(n-1)}
{n(n-1)+(n/2+1)(n/2)}
=
\frac{8(n-1)}{5n-2}.
\]
Thus,
\[
m(F_{n/2}(G))
\leq
\frac{8(n-1)}{5n-2}
\binom{n-2}{n/2-1}m(G)
=
\left(\frac85+o(1)\right)
\binom{n-2}{n/2-1}m(G);
\]
which proves $c)$.
\end{proof}

Note that the multiplicative constant gap between $m(F_k(G))$ and  (\ref{eq:lower}) proven in Theorem~\ref{thm:K5q} holds when $k=\Theta(n)$;
 Theorem~\ref{thm:upper_bound} shows that the lower bound of (\ref{eq:lower}) is asymptotically tight when $k=o(n)$.
\section{Bounds for some classes of parent graphs }\label{sec:parents}

\subsection{Complete Graph}

Throughout this subsection assume that $G$ is the complete graph $K_n$ on $n$ vertices, and that its vertex set is equal to $[n]$.
Let $(X',Y')$ be a cut of $F_k(K_n)$, $A \in X'$ and $B \in Y'$. Note that
$AB \in E(X',Y')$ if and only if $A \cap B \in \binom{[n]}{k-1}$.
Thus, if  for every  $(k-1)$-subset $S$, the $k$-subsets of $[n]$ containing $S$ are divided as
evenly as possible between $X'$ and $Y'$, then $(X',Y')$ is maximum.
 We use this observation to prove the following
upper bound on $m(F_k(K_n))$.

\begin{theorem}\label{mxcut_johnson}
    \begin{align}
      m(F_k(K_n)) & \leq\left\lfloor\frac{(n-k+1)^2}{4}\right\rfloor \binom{n}{k-1} \label{eq:upKn}\\[8pt]
      & \le \left(1+\frac{2}{n}\right)\binom{n-2}{k-1}m(K_n). \label{eq:knasym}
    \end{align}
Where equality in (\ref{eq:upKn}) is attained if and only if there is a cut $(X',Y')$ of $F_k(K_n)$  such that
for every $S\in\binom{[n]}{k-1}$,
\[\{A \in X':S \subset A\} \textrm{ and } \{B \in Y':S \subset B\}\]
are of sizes
$\left\lceil\frac{n-k+1}{2}\right\rceil$ and
$\left\lfloor\frac{n-k+1}{2}\right\rfloor$.
\end{theorem}
\begin{proof}
From the observation preceding this theorem, we have that
 \[m(F_k(K_n)) \le \left\lceil\frac{n-k+1}{2}\right\rceil\left\lfloor\frac{n-k+1}{2}\right\rfloor \binom{n}{k-1} =\left\lfloor\frac{(n-k+1)^2}{4}\right\rfloor \binom{n}{k-1},\]
 and the condition for achieving equality.

 It remains to show (\ref{eq:knasym}). Let $N:=n-k+1$. Since
 \[
m(K_n)=\left\lfloor\frac{n^2}{4}\right\rfloor\geq\frac{n^2-1}{4}
\textrm{ and }
\binom{n}{k-1}=\frac{n(n-1)}{N(N-1)}\binom{n-2}{k-1},
\]
we have
\begin{align*}
\binom{n}{k-1}\left\lfloor\frac{N^2}{4}\right\rfloor
&=
\frac{n(n-1)}{N(N-1)}\,\binom{n-2}{k-1}
\left\lfloor\frac{N^2}{4}\right\rfloor\\
&\leq
\frac{n(n-1)}{N(N-1)}\,\binom{n-2}{k-1}\cdot
\frac{N^2}{n^2-1}\cdot\frac{n^2-1}{4}\\
&\leq
\frac{n(n-1)}{N(N-1)}\cdot\frac{N^2}{n^2-1}\,
\binom{n-2}{k-1}\,m(K_n)\\
&=
\frac{N}{N-1}\cdot\frac{n}{n+1}\,
\binom{n-2}{k-1}\,m(K_n).
\end{align*}
Finally, as we are assuming $k\le n/2$,
\[
\frac{N}{N-1}\cdot\frac{n}{n+1}
\leq
\frac{N}{N-1}
=
1+\frac{1}{N-1}
=
1+\frac{1}{n-k}
\leq
1+\frac{2}{n}.
\]
The result follows.
\end{proof}

For the following lower bounds we use Walecki's construction~\cite{Lucas1892,Alspach2008}.
He showed that for odd $n$, the complete graph $K_n$ decomposes into $(n-1)/2$
Hamiltonian cycles, and for even $n$ it decomposes into $(n-2)/2$
Hamiltonian cycles together with a perfect matching.

\begin{prop}\label{prop:exact_small_k}
\[
m(F_2(K_n))
=
\begin{cases}
n\left\lfloor\dfrac{(n-1)^2}{4}\right\rfloor,
& n\not\equiv 3 \pmod 4,\\[4ex]
n\left\lfloor\dfrac{(n-1)^2}{4}\right\rfloor-1,
& n\equiv 3\pmod 4.
\end{cases}
\]
\end{prop}

\begin{proof}

A cut of $F_2(K_n)$ is a partition of $E(K_n)$ into two
classes $X'$ and $Y'$. By
Theorem~\ref{mxcut_johnson},
\[
m(F_2(K_n))
\leq
n
\left\lfloor\frac{(n-1)^2}{4}\right\rfloor,
\]
and the upper bound is attained precisely
when, at every vertex $v$ of $K_n$, the number of edges of $X'$ and $Y'$ incident to $v$ are equal to $\left\lceil\frac{n-1}{2}\right\rceil$ and
$\left\lfloor\frac{n-1}{2}\right\rfloor$.

Suppose that $n\equiv3\pmod4$. Thus, $n$ is odd and $n-1$ is even. Equality in
the upper bound requires that every vertex is incident to $(n-1)/2$ edges in each of $X'$ and $Y'$. But each of $X'$ and $Y'$
would be the edge set of an $(n-1)/2$-regular spanning subgraph,
whose sum of degrees would be equal to $n\cdot\frac{n-1}{2}$, which is odd---a contradiction.
Therefore,
\[
m(F_2(K_n))
\leq
n
\left\lfloor\frac{(n-1)^2}{4}\right\rfloor-1,
\]
in this case.
We now construct the cut $(X',Y')$  satisfying the upper bounds in each case.

Suppose that $n \equiv 1 \pmod 4$. In this case $(n-1)/4$ is an integer. Assign the edges of half of the $(n-1)/2$ Hamiltonian cycles in Walecki's construction to $X'$, and the edges of the other half to $Y'$. Every vertex of $K_n$ is incident to $\frac{n-1}{2}$ edges of each of $X'$ and $Y'$.

Suppose  that $n$
is even.  Assign the edges of $\lceil (n-2)/4\rceil$ of the Hamiltonian cycles in Walecki's construction  to $X'$, and the remaining together with the perfect matching
to $Y'$. Every vertex of $K_n$ is incident to
$2\lceil (n-2)/4\rceil$ edges in $X'$ and
$2\lfloor (n-2)/4\rfloor+1$ edges in $Y'$; these two numbers differ by
exactly one; thus, one of them is equal to $\left\lceil\frac{n-1}{2}\right\rceil$
and the other to $\left\lfloor\frac{n-1}{2}\right\rfloor$.

Suppose that $n\equiv3\pmod4$. Thus, $(n-3)/4$ is an integer. Let $X'$ be
the union of the edges of $(n-3)/4$ of the Hamiltonian cycles in Walecki's construction together with a near-perfect matching
of one of the remaining cycles. Let $Y':=E(K_n)\setminus X'$.
Thus, every vertex but one is incident to $\frac{n-3}{2}+1=\frac{n-1}{2}$ edges of $X'$, and to $\frac{n-1}{2}$ edges of $Y'$. The remaining vertex is incident to
$\frac{n-3}{2}$ edges of $X'$ and $\frac{n+1}{2}$ edges of $Y'$.
By (\ref{eq:cut_grouping}), $(X',Y')$ has weight equal to
\[
(n-1)\left(\frac{n-1}{2}\right)^2
+\frac{n-3}{2}\cdot\frac{n+1}{2}
=
n\left(\frac{n-1}{2}\right)^2-1=n
\left\lfloor\frac{(n-1)^2}{4}\right\rfloor-1.
\]
\end{proof}

 A \emph{$t$-$(v,k,\lambda)$ design} is a family of $k$-subsets, called \emph{blocks}, of a $v$-set in which
every $t$-subset is contained in exactly $\lambda$ blocks. When  the family consists
of all $k$-subsets of $[n]$ we call it the \emph{complete design}; this is a
design with $\lambda=\binom{n-t}{k-t}$.
A partition of $\binom{[n]}{k}$ into two sets, in which each set is a $t$-$\left (n,k,\frac{1}{2}\binom{n-t}{k-t} \right )$ design, is called
a \emph{halving of the complete design}; see~\cite{hartman1987halving}.
When $(X',Y')$ attains the upper bound of Theorem~\ref{mxcut_johnson}
and $n-(k-1)$ is even,
each of $X'$ and $Y'$, is a $(k-1)-(n,k,(n-k+1)/2)$-design;
and $(X',Y')$ is a halving of the complete design.
We use this relationship with block designs to prove the exact value of
$m(F_3(K_n))$.

\begin{prop}\label{thm:F3}
\[
    m(F_3(K_n))=
    \begin{cases}
        \displaystyle\binom{n}{2}\frac{(n-2)^2}{4}-\frac{n}{2}
            & \textrm{if } n\equiv 0\pmod 4,\\[4ex]
        \displaystyle\binom{n}{2}\left\lfloor\frac{(n-2)^2}{4}\right\rfloor
            & \textrm{otherwise.}
    \end{cases}
\]
%
%
\end{prop}
\begin{proof}
We first prove the upper bound.
By Theorem~\ref{mxcut_johnson} with $k=3$,
\begin{equation} \label{eq:k3}
    m(F_3(K_n))\leq\binom{n}{2}\left\lfloor\frac{(n-2)^2}{4}\right\rfloor,
\end{equation}
with equality if and only if some cut $(X',Y')$ of $F_3(K_n)$ is such that for every $S\in\binom{[n]}{2}$, the number of triples containing $S$, in $X'$ and $Y'$ are equal to  $\left\lceil\frac{n-2}{2}\right\rceil$ and $\left\lfloor\frac{n-2}{2}\right\rfloor$.

Suppose that $n\equiv0\pmod4$, and let $(X',Y')$ be a cut of $F_3(K_n)$.
Let $S\in\binom{[n]}{2}$, and let $d$ be the number of triples of $X'$ containing $S$; each of these triples
is adjacent each of the triples of $Y'$ that contain $S$.
Thererefore, the number of edges $AB\in E(X',Y')$,  with $A\cap B=S$, is equal to $d(n-2-d)$; this is equal to $\frac{(n-2)^2}{4}$ if $d=\frac{n-2}{2}$, and
at most $\frac{(n-2)^2}{4}-1$ otherwise.
Let $v\in[n]$, and consider the graph on $[n]\setminus\{v\}$ in which $uw$ is an edge
whenever $\{v\} \cup \{u,w\}\in X'$. Note that the degree of a vertex
$u$ in this graph is the number of
triples of $X'$ containing $\{v,u\}$. Since $n-1$ and $\frac{n-2}{2}$ are both odd,
this graph is not $\frac{n-2}{2}$-regular. Thus, for every $v\in[n]$ there exists a
$u\in[n]\setminus\{v\}$ such that the number of triples of $X'$ containing $\{v,u\}$
is not equal to $\frac{n-2}{2}$. These pairs cover $[n]$, so there are at least
$\frac{n}{2}$ of them. Therefore,
\begin{equation}\label{eq:k3_0}
    m(F_3(K_n)) \leq\binom{n}{2}\frac{(n-2)^2}{4}-\frac{n}{2},
\end{equation}
in this case.

We now prove the lower bound.
Dehon~\cite{dehon} showed that $\binom{[r]}{3}$ can be
halved if and only if $r\equiv2\pmod4$.
Note that when $n$ is odd, equality in (\ref{eq:k3}) requires that every $S\in\binom{[n]}{2}$
is in $\frac{n-1}{2}$ triples on one of $X'$ or $Y'$, and in  $\frac{n-3}{2}$ of the
other.
\begin{itemize}
 \item  Suppose that $n\equiv1\pmod4$.

 Thus, $\binom{[n+1]}{3}$ can be halved. Let $X''$ be one of
the two halves, and set $X':=X''\cap\binom{[n]}{3}$ and $Y'=\binom{[n]}{3}\setminus X''$. Every $S\in\binom{[n]}{2}$ lies in exactly $\frac{n-1}{2}$ triples of $X''$,
of which at most one, namely $S\cup\{n+1\}$, contains the point $n+1$.
Therefore $S$ lies in $\frac{n-1}{2}$ or $\frac{n-3}{2}$ triples of
$X'$, as required.

\item Suppose that $n\equiv3\pmod4$.

Thus, $\binom{[n-1]}{3}$ can be halved. Let $X''$ be one of
the two halves. Note that every $S \in \binom{[n-1]}{2}$
lies in exactly $\frac{n-3}{2}$ elements of $X''$.
Let $R$ be the subgraph of $K_{n-1}$ obtained by  taking $\frac{n-3}{4}$ of the
Hamiltonian cycles in Walecki's construction together with the perfect matching. Let
\[
X':=X''\cup\bigl\{\{x,y,n\}:xy\in E(R)\bigr\} \textrm{ and }
Y':=\binom{[n]}{3}\setminus X'.
\]
A given $\{x,y\}\in \binom{[n-1]}{2}$ lies in $\frac{n-3}{2}$ or $\frac{n-1}{2}$
elements of $X'$, depending on whether $xy$ lies in $R$; and
$\{x,n\} \in \binom{[n]}{2}$ lies in as many elements of $X'$ as $x$ has neighbors in
$R$, that is, in $\frac{n-1}{2}$; so $(X',Y')$ are  as required.

\item Suppose that $n\equiv2\pmod4$.

Equality in (\ref{eq:k3}) requires
that every $S\in\binom{[n]}{2}$ is contained in exactly $\frac{n-2}{2}$ triples of each of $X'$
and $Y'$. This is a halving of $\binom{[n]}{3}$ into two $2$-$(n,3,\frac{n-2}{2})$ designs. Such a halving exists by Dehon's result.

\item Suppose that $n\equiv0\pmod4$.

Thus, $\binom{[n-2]}{3}$ can be halved. Let $X''$ be one of
the two halves. Note that every $S \in \binom{[n-2]}{2}$
lies in exactly $\frac{n-4}{2}$ elements of $X''$.
Walecki's construction decomposes $K_{n-2}$ into $\frac{n-4}{2}$ Hamiltonian cycles
and a perfect matching $M$. Let $M_0$ be a set of $\frac{n-4}{4}+1$ edges of $M$, and
let $M_1$ be the set of the $\frac{n-4}{4}$ remaining edges of $M$.
Let $R_0$ be the subgraph of $K_{n-2}$ obtained by taking $\frac{n-4}{4}$ of the
Hamiltonian cycles together with $M_0$, and let $R_1$ be the subgraph obtained by taking
the remaining $\frac{n-4}{4}$ Hamiltonian cycles together with $M_0$. Let
\begin{align*}
X':=X'' & \cup\bigl\{\{x,y,n-1\}:xy\in E(R_0)\bigr\}
          \cup\bigl\{\{x,y,n\}:xy\in E(R_1)\bigr\}\\
        & \cup\bigl\{\{x,n-1,n\}:x \textrm{ is an endpoint of an edge of } M_1\bigr\}
\end{align*}

and \[Y':=\binom{[n]}{3}\setminus X'.\]
Let $\{x,y\}\in \binom{[n-2]}{2}$.
If $xy$ is an edge of one of the $\frac{n-4}{2}$ Hamiltonian cycles in Walecki's construction, then $xy$ is an edge of exactly one of $R_0$ and $R_1$; thus, $\{x,y\}$ is
contained in exactly $\frac{n-4}{2}+1=\frac{n-2}{2}$ elements of $X'$.
If $xy \in M_0$, then $xy$ is an edge of both $R_0$ and $R_1$; thus, $\{x,y\}$ is contained
in exactly $\frac{n-4}{2}+2=\frac{n-2}{2}+1$ elements of $X'$.
If $xy \in M_1$, then $xy$ is an edge of neither $R_0$ nor $R_1$; thus, $\{x,y\}$ is contained
in exactly $\frac{n-4}{2}=\frac{n-2}{2}-1$ elements of $X'$.

Let $x \in [n-2]$. Note that $x$ is an endpoint of exactly one edge of $M$.
If this edge is in $M_0$, then $x$ has degree
$2\cdot\frac{n-4}{4}+1=\frac{n-2}{2}$ in each of $R_0$ and $R_1$, and
$\{x,n-1,n\}\notin X'$. If this edge is in $M_1$, then $x$ has degree
$2\cdot\frac{n-4}{4}=\frac{n-4}{2}$ in each of $R_0$ and $R_1$, and
$\{x,n-1,n\}\in X'$. In both cases, each of $\{x,n-1\}$ and $\{x,n\}$ is contained in
exactly $\frac{n-2}{2}$ elements of $X'$.

Finally, $\{n-1,n\}$ is contained in exactly $2|M_1|=\frac{n-4}{2}$ elements of $X'$.

Thus, every $S\in\binom{[n]}{2}$ is contained in exactly $\frac{n-2}{2}$ elements of $X'$,
except for  $\{n-1,n\}$ and the $\frac{n-2}{2}$ edges of $M$. Each of these
$\frac{n}{2}$ pairs is contained in exactly $\frac{n-2}{2}+1$ or $\frac{n-2}{2}-1$
elements of $X'$. Therefore,
\begin{align*}
 |E(X',Y')|&=\left (\binom{n}{2}-\frac{n}{2} \right ) \left(\frac{n-2}{2}\cdot \frac{n-2}{2} \right)+\frac{n}{2}\left ( \frac{n-2}{2}+1\right) \left ( \frac{n-2}{2}-1\right) \\
 &=\binom{n}{2}\frac{(n-2)^2}{4}-\frac{n}{2}
\end{align*}

and equality holds in (\ref{eq:k3_0}).

\end{itemize}
\end{proof}

\subsection{Triangular and Kagome graphs}\label{sec:lattices}

Geometrically frustrated spin systems are most often studied on regular
two-dimensional lattices whose interaction graph is not bipartite. Two of the
most prominent examples are the \emph{triangular} and the \emph{Kagome}
lattices, and the physical interest lies in finite pieces of them. Let
\[
\mathbf{e}_1:=(1,0)
\quad\text{and}\quad
\mathbf{e}_2:=\Bigl(\tfrac12,\tfrac{\sqrt3}{2}\Bigr).
\]

\begin{definition}[Triangular lattice]
The \emph{triangular lattice} $\mathbb{T}$ is the infinite graph with vertex
set $\{\,a\,\mathbf{e}_1+b\,\mathbf{e}_2 : a,b\in\mathbb{Z}\,\}$ in which two
vertices are adjacent if and only if they lie at Euclidean distance equal to $1$.
\end{definition}

\begin{figure}[ht]
\centering
\begin{subfigure}[b]{0.48\textwidth}
\centering
\includegraphics[width=\textwidth]{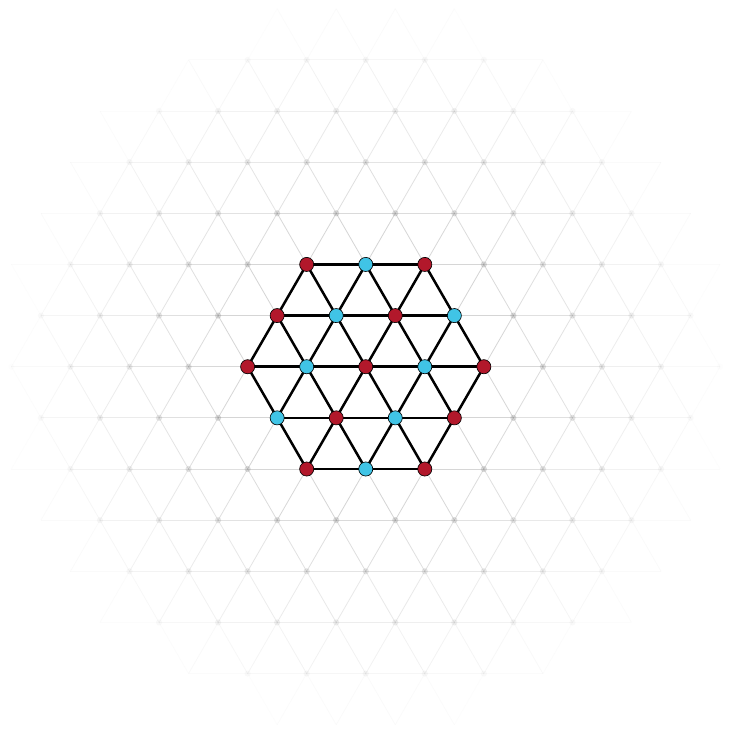}
\caption{A triangular-lattice graph.}
\label{fig:triangular}
\end{subfigure}
\hfill
\begin{subfigure}[b]{0.48\textwidth}
\centering
\includegraphics[width=\textwidth]{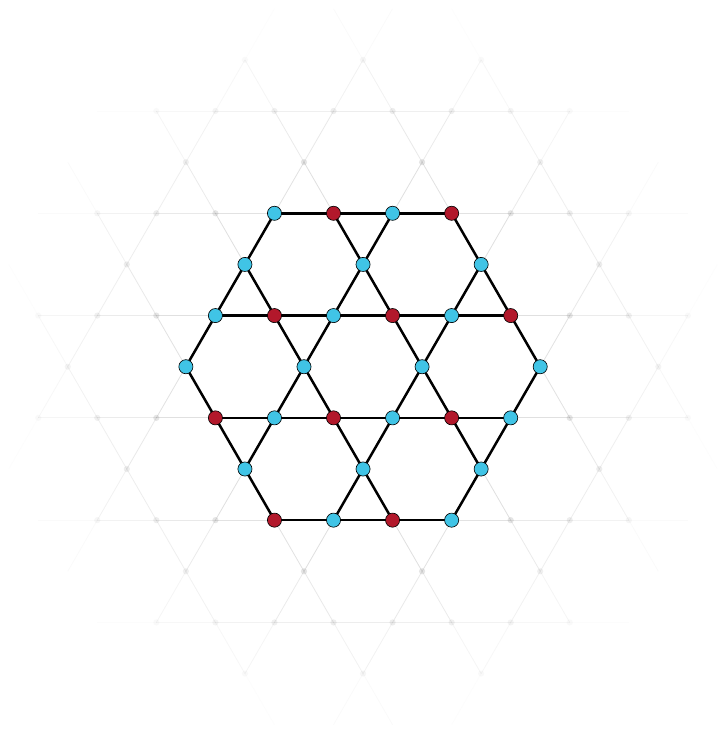}
\caption{A Kagome graph.}
\label{fig:kagome}
\end{subfigure}
\caption{Finite induced subgraphs of the triangular and Kagome lattices; the
remainder of each lattice is shown faded in gray. A vertex
$a\,\mathbf{e}_1+b\,\mathbf{e}_2$ is colored red if $a$ is even and blue if
$a$ is odd, so that the two colors alternate along every horizontal row of $\mathbb{T}$. }
\label{fig:lattices}
\end{figure}

\begin{definition}[Kagome lattice]
The \emph{Kagome lattice} $\mathbb{K}$ is the infinite graph with vertex set
$\{\,2a\,\mathbf{e}_1+2b\,\mathbf{e}_2+\mathbf{c} : a,b\in\mathbb{Z},\
\mathbf{c}\in\{\mathbf{e}_1,\mathbf{e}_2,\mathbf{e}_1+\mathbf{e}_2\}\,\}$ in
which two vertices are adjacent if and only if they lie at Euclidean distance equal to
$1$. Equivalently, $\mathbb{K}$ is obtained from
$\mathbb{T}$ by deleting the sublattice
$\{\,2a\,\mathbf{e}_1+2b\,\mathbf{e}_2 : a,b\in\mathbb{Z}\,\}$.
\end{definition}
A \emph{triangular-lattice graph} is a graph isomorphic to a finite induced
subgraph of $\mathbb{T}$, and a \emph{Kagome graph} is a graph isomorphic to a
finite induced subgraph of $\mathbb{K}$. See Figure~\ref{fig:lattices}.

It is easy to show (see \cite{Token}) that $F_k(G)$ is bipartite if and only if $G$ is bipartite. The maximum cut of a bipartite graph is given by the bipartition of its vertices. Triangular lattice and Kagome graphs may have triangles; thus, in general, they are not bipartite.
However, the triangles of these lattices have a
structure that allows us to find near-tight bounds on $m(F_k(G))$ in many instances for their induced finite subgraphs.

We consider a more general setting. Let $\mathcal{T}$ be a
set of pairwise edge disjoint triangles of  $G$. Since at most two edges of
every triangle lie in a given cut set, and the triangles of $\mathcal{T}$ share
no edges, we have that
\begin{equation}\label{eq:triangles_upper_bound}
    m(G) \le |E(G)|-|\mathcal{T}|.
\end{equation}
In some instances we can bound $m(G)$ from below in terms of $|\mathcal{T}|$ as
well.
\begin{prop}\label{prop:hetero}
 Suppose that there is a two-coloring of the vertices of $G$ such that
 every triangle in $\mathcal{T}$ is bichromatic. Then
 \[m(G) \ge 2|\mathcal{T}|.\]
 In particular, if $E(\mathcal{T})=E(G)$, then
 \[m(G) = 2|\mathcal{T}|.\]
\end{prop}
\begin{proof}
Let $(X,Y)$ be the partition of $V(G)$ induced by the color classes. Since every
triangle of $\mathcal{T}$ is bichromatic, each triangle of $\mathcal{T}$ contributes two edges to
$E(X,Y)$, giving the
first inequality. If  $E(\mathcal{T})=E(G)$ then
$|E(G)|=3|\mathcal{T}|$, and (\ref{eq:triangles_upper_bound}) gives the reverse
inequality.
\end{proof}

We color the vertices of  $\mathbb{T}$ and  $\mathbb{K}$ as follows.
Color
$a\,\mathbf{e}_1+b\,\mathbf{e}_2$ red if $a$ is even and blue if $a$ is odd. Note
that the colors alternate along each horizontal row of $\mathbb{T}$. Thus, every
triangle of $\mathbb{T}$ is bichromatic. Passing from $\mathbb{T}$ to $\mathbb{K}$
deletes precisely the red vertices of the rows with $b$ even, and leaves the
rows with $b$ odd untouched; since every triangle of $\mathbb{K}$ has an edge
contained in one of these latter rows, its endpoints again receive different
colors. See Figure~\ref{fig:lattices}.

In \cite{Token} it is shown that there are exactly two ways of
 constructing a triangle $T'$ in $F_k(G)$, either
    \begin{enumerate}
        \item by moving one token in a triangle $T$ of $G$, while fixing the remaining $k-1$ tokens in a set of vertices $S \subset V(G) \setminus V(T)$, or
        \item  by moving two tokens in a triangle $T$ of $G$, while fixing the remaining $k-2$ tokens in a set of vertices $S \subset V(G) \setminus V(T)$.
    \end{enumerate}
\noindent We call $T$ the \emph{parent} of $T'$. Note that $T$ is parent of  exactly $\binom{n-3}{k-1}+\binom{n-3}{k-2}= \binom{n-2}{k-1}$
triangles in $F_k(G)$.
Let $\mathcal{T'}$ be the triangles of $F_k(G)$ whose parents are in $\mathcal{T}$. We have that
\begin{equation}\label{eq:triangle_lift}
 \text{$\mathcal{T'}$ is a set of $|\mathcal{T}|\binom{n-2}{k-1}$  pairwise edge disjoint triangles of $F_k(G)$.}
\end{equation}

A red-blue coloring of the vertices of $G$ lifts to a red-blue coloring of the vertices of $F_k(G)$ as follows.
Color a token configuration $A$ red if $A$ has an
even number of red vertices, and blue otherwise. A triangle of $F_k(G)$ is bichromatic
under this coloring if and only if its parent triangle is bichromatic in
$G$.
In particular, lifting the colorings of $\mathbb{T}$ and $\mathbb{K}$ given
above, we obtain, for every triangular or Kagome graph $G$, a
coloring of $V(F_k(G))$ in which every triangle of $F_k(G)$ is
bichromatic.

A triangular-lattice or a Kagome graph $G$ is a \emph{patch} if its outer
face is bounded by a cycle $C$, and  $G$ contains every lattice
vertex enclosed by $C$.
In particular, the \emph{triangular patch} of side $m$ is the triangular lattice  patch graph, $\Delta_m$, whose outer face is an equilateral triangle of side length equal to $m$.
It can be shown that
\[
|V(\Delta_m)|=\frac{(m+1)(m+2)}{2}\quad \text{and} \quad |E(\Delta_m)|= \frac{3m(m+1)}{2}.
\]

\begin{prop}\label{prop:Tm}
    Let $n:=|V(\Delta_m)|$. Then
    \[
    m(F_k(\Delta_m))=m(m+1)\binom{n-2}{k-1}.
    \]
\end{prop}
\begin{proof}
Let $\mathcal{T}$ be the set of the $\frac{m(m+1)}{2}$ upward unit triangles of
$\Delta_m$. They are pairwise edge disjoint and together use all of the
$\frac{3m(m+1)}{2}$ edges of $\Delta_m$, so $E(\mathcal{T})=E(\Delta_m)$.
Lifting them as in (\ref{eq:triangle_lift}) gives a set $\mathcal{T}'$ of
$\frac{m(m+1)}{2}\binom{n-2}{k-1}$ pairwise edge disjoint triangles of
$F_k(\Delta_m)$, which is a third of
\[
|E(F_k(\Delta_m))|=|E(\Delta_m)|\binom{n-2}{k-1}
=\frac{3m(m+1)}{2}\binom{n-2}{k-1};
\]
so $E(\mathcal{T}')=E(F_k(\Delta_m))$. Finally, every triangle of
$F_k(\Delta_m)$ is bichromatic in the lifted coloring of $\mathbb{T}$. Therefore,
by Proposition~\ref{prop:hetero}, we have
\[
m(F_k(\Delta_m))=2\,|\mathcal{T}'|=m(m+1)\binom{n-2}{k-1}.
\]
\end{proof}

\begin{theorem}\label{thm:polygon}
 Suppose that $G$ is a patch (triangular or Kagome). Let   $L$ be the length of its outer face $C$. Then
\[
m(F_k(G))
\;\leq\;
\left(1+\frac{2L}{n}\right)
\binom{n-2}{k-1}m(G).
\]
Thus, if $L=o(n)$, then
\[
m(F_k(G))=(1+o(1))\binom{n-2}{k-1}m(G).
\]
\end{theorem}

\begin{proof}

Let $\mathbb{L}$ be  the lattice, $\mathbb{T}$ or $\mathbb{K}$, of which $G$
is a finite induced subgraph; and let $R$ be the closed region bounded by $C$.
There are $L$ vertices of $G$ lying in the boundary of $R$ and $n-L$
vertices in its interior.
If $\mathbb{L}=\mathbb{K}$, then let $\mathcal{T}$ be the set of triangles of $G$;
if $\mathbb{L}=\mathbb{T}$, then let $\mathcal{T}$ be the set of upward triangles of $G$.
Note that:
if $\mathbb{L}=\mathbb{T}$, then every interior vertex of $G$ is
the leftmost vertex of exactly one triangle in $\mathcal{T}$;
 if $\mathbb{L}=\mathbb{K}$, then every interior vertex of $G$ is
in exactly two triangles in $\mathcal{T}$, and
every edge of $G$, in the interior of $R$, is contained in exactly one triangle of $\mathcal{T}$.
In both cases we have that
\[
|\mathcal{T}| \geq \frac{2}{3}(n-L) \text{ and } |E(G)| \leq 3|\mathcal{T}|+L.
\]

We assume that $n\geq 2L$, since otherwise $2L/n>1$ and the upper bound  holds trivially, given that
\[
m(F_k(G))\leq |E(F_k(G))|=\binom{n-2}{k-1}|E(G)|\leq 2\binom{n-2}{k-1}m(G);
\]
 thus, $|\mathcal{T}|\geq n/3$ for both lattices.

 Every triangle of $\mathcal{T}$ is bichromatic for the coloring described
above, so by Proposition~\ref{prop:hetero} we have that
$m(G)\geq 2|\mathcal{T}|\geq 2n/3$.
By (\ref{eq:triangle_lift}),  $F_k(G)$ has a set of
$|\mathcal{T}|\binom{n-2}{k-1}$ pairwise edge disjoint triangles. Thus,
by (\ref{eq:triangles_upper_bound}) we have that
\[
\begin{aligned}
m(F_k(G))
&\leq
\binom{n-2}{k-1}\bigl(|E(G)|-|\mathcal{T}|\bigr)
\leq
\binom{n-2}{k-1}\bigl(2|\mathcal{T}|+L\bigr)\\
&\leq
\binom{n-2}{k-1}\bigl(m(G)+L\bigr)
=
\left(1+\frac{L}{m(G)}\right)\binom{n-2}{k-1}m(G).
\end{aligned}
\]
Finally, $m(G)\geq 2n/3$ gives
\[
\frac{L}{m(G)}\leq\frac{3L}{2n}\leq\frac{2L}{n}.
\]
The result follows.
\end{proof}

\subsection{Random graphs}\label{sec:random}

The upper bound (\ref{eq:cut_grouping_ub})  of $\sum m(G-S)$ on $m(F_k(G))$ is obtained
by considering for every  $S \in \binom{V(G)}{k-1}$ the projections $(X_S,Y_S)$ of a given cut $(X',Y')$ of  $F_k(G)$. Intuitively, the bound
should be near tight when the token configurations containing these $k-1$ sets are well spread among  $X'$  and $Y'$. This suggests that the parity lift should give a near maximum cut in the case of random graphs.

Throughout this section $ G(n,p)$ denotes the
Erd\H{o}s--R\'enyi model of random graphs; in this model there are $n$ vertices and every pair is made adjacent
independently with probability equal to $p$. We write $G \sim G(n,p)$ for  a random graph from this model.
We allow $p=p(n)$ to depend on $n$. In this context \emph{with
high probability (w.h.p.)} means with probability tending to $1$ as $n\to\infty$.

 We call
\[
\sigma(G):=m(G)-\frac{|E(G)|}{2}\geq 0
\]
the \emph{surplus} of $G$.
\begin{prop}\label{prop:surplus}
\[
m(F_k(G))
\;\leq\;
\left(1+\frac{3\,\sigma(G)}{m(G)}\right)\binom{n-2}{k-1}m(G).
\]
\end{prop}

\begin{proof}
Substituting $w(G) = |E(G)|=2m(G)-2\sigma(G)$ into the second expression inside the
minimum of (\ref{eq:minbound}) gives
\[
\begin{aligned}
m(F_k(G))
&\leq
\binom{n}{k-1}m(G)
-\left(\binom{n}{k-1}-\binom{n-2}{k-1}\right)\bigl(m(G)-\sigma(G)\bigr)\\
&=
\binom{n-2}{k-1}m(G)
+\left(\binom{n}{k-1}-\binom{n-2}{k-1}\right)\sigma(G).
\end{aligned}
\]
Since $n-k\geq n/2$, we have that
\[
\frac{\binom{n}{k-1}}{\binom{n-2}{k-1}}
=
\frac{n(n-1)}{(n-k+1)(n-k)}
<
\frac{n^2}{(n/2)^2}
=
4.
\]
Thus $\binom{n}{k-1}-\binom{n-2}{k-1}<3\binom{n-2}{k-1}$, and the result follows.
\end{proof}

\begin{theorem}\label{thm:random}
Let $G\sim G(n,p)$ with $pn\to\infty$. Then, w.h.p,
\[
m(F_k(G))
=
\left(1+O\!\left(\frac{1}{\sqrt{pn}}\right)\right)
\binom{n-2}{k-1}m(G).
\]
\end{theorem}

\begin{proof}

Let $S\subseteq[n]$. Note that
\[|E_G(S,V(G) \setminus S)|=\sum_{\substack{u \in S \\ v \in V(G)\setminus S}} X_{(u,v)},\]
where the   $X_{(u,v)}$ are independent Bernoulli random variables with parameter $p$. For every
 $u \in S$ and $v \in V(G) \setminus S$, let $Y_{(u,v)}:=X_{(u,v)}-p$; this is an
 independent random variable of mean equal to zero, $|Y_{(u,v)}|\leq 1$, and variance equal to
 $p(1-p)$. Therefore, the total variance is equal to
 \[
|S|(n-|S|)\,p(1-p)\leq \frac{p(1-p)n^{2}}{4}.
\]

By Bernstein's inequality for sums of independent, zero mean and bounded random
variables, we have that
\[
\mathbb P\!\left(\sum_{u\in S}\sum_{v\in V(G)\setminus S}Y_{(u,v)}\geq t\right)
\leq
\exp\!\left(-\frac{t^{2}}{p(1-p)n^{2}/2+2t/3}\right),
\]
for every $t\geq 0$.
Since
\[
\sum_{u\in S}\sum_{v\in V(G)\setminus S}Y_{(u,v)}
=
\sum_{u\in S}\sum_{v\in V(G)\setminus S}X_{(u,v)}-p\,|S|(n-|S|)
=
|E_G(S,V(G)\setminus S)|-p\,|S|(n-|S|),
\]
this may be rewritten as
\begin{equation}\label{eq:bernstein}
\mathbb P\bigl(|E_G(S,V(G)\setminus S)|-p\,|S|(n-|S|)\geq t\bigr)
\leq
\exp\!\left(-\frac{t^{2}}{p(1-p)n^{2}/2+2t/3}\right).
\end{equation}
Let  $t=2\bigl(\sqrt{p}\,n^{3/2}+n\bigr)$.  We have that
\[
\begin{aligned}
t^{2}
&=
4\left(pn^{3}+2\sqrt{p}\,n^{5/2}+n^{2}\right) \\
& >
pn^{3}+\frac{8}{3}\sqrt{p}\,n^{5/2}+\frac{8}{3}n^{2}\\
&\geq
p(1-p)n^{3}+\frac{8}{3}\sqrt{p}\,n^{5/2}+\frac{8}{3}n^{2}\\
&=
2n\left(\frac{p(1-p)n^{2}}{2}+\frac{2t}{3}\right).
\end{aligned}
\]

By (\ref{eq:bernstein}), for a given $S$, the probability that
$|E_G(S,V(G)\setminus S)|$ exceeds $p\,|S|(n-|S|)+t$ is at most $e^{-2n}$. Since there
are at most $2^{n}$ sets
$S$, and $2^{n}e^{-2n}\to 0$, the union bound shows that w.h.p, every $S$ satisfies that
$|E_G(S,V(G)\setminus S)|\leq p\,|S|(n-|S|)+t$. Since $p\,|S|(n-|S|)\leq pn^{2}/4$ and $pn \to \infty$, this
gives
\begin{equation}\label{eq:cut_ub}
m(G)=\max_{S\subseteq[n]}|E_G(S,V(G)\setminus S)|
\leq
\frac{pn^{2}}{4}+2\left(\sqrt{p}\,n^{3/2}+n\right)
=\frac{pn^{2}}{4}+O\!\left(\sqrt{p}\,n^{3/2}\right)
\end{equation}
w.h.p.

By applying Bernstein's inequality in a similar way with the hypothesis $pn \to \infty$, we obtain that w.h.p
\begin{equation}\label{eq:edges_lb}
\frac{|E(G)|}{2}
\geq
\frac{pn(n-1)}{4}-\frac{\sqrt{p}\,n\log n}{2}
\geq
\frac{pn^{2}}{4}-O\!\left(\sqrt{p}\,n^{3/2}\right).
\end{equation}
Therefore,
\[
\sigma(G)=m(G)-\frac{|E(G)|}{2}=O\!\left(\sqrt{p}\,n^{3/2}\right).
\]

Since $m(G)\geq|E(G)|/2$ and $pn \to \infty$, we have that
$m(G)=\Omega\bigl(pn^{2}\bigr)$ and
\[
\frac{\sigma(G)}{m(G)}
=
O\!\left(\frac{\sqrt{p}\,n^{3/2}}{pn^{2}}\right)
=
O\!\left(\frac{1}{\sqrt{pn}}\right);
\]
the result follows from Proposition~\ref{prop:surplus}.
\end{proof}

\section{Algorithmic complexity}\label{sec:complexity}

We now consider the algorithmic complexity of computing $m(F_k(G))$.
We consider the following three scenarios depending on the problem's input.

\begin{itemize}
\item \emph{Token generated.} The input is the graph $F_k(G)$, each of whose
vertices is labeled by the token configuration it represents. Recovering $G$
from such an input is straightforward; letting $N:=\binom{n}{k}$, the input has size $O(n^2N)$.

\item \emph{Token isomorphic.} The input is a graph $H$ isomorphic to a token
graph $F_k(G)$, but the isomorphism is not given. We do not know whether such
an isomorphism can be found in time polynomial in $|H|$. In the previous problem
having the labeling at hand may allow for faster algorithms.

\item \emph{Token parent.} The input is a graph $G$ and an integer $k$, and the
questions asked are about $F_k(G)$. In some instances it is possible to solve the algorithmic problem
without constructing $F_k(G)$.
\end{itemize}

Although \textsc{Max-Cut} is in NP when the graph is given explicitly,
it is not clear this is the case in the token parent scenario;
since a certificate
might have superpolynomial size in $n$ if, say, $k=\omega(1)$.
We consider the \textsc{Max-Cut} problem in the token generated scenario.

\begin{problem}[\textsc{Token Max-Cut}]\label{prob:tmc}
Given a token graph $F_k(G)$, with its vertices labeled by the corresponding
token configurations and its edges given explicitly, together with a weight
function $w:E(F_k(G))\to \mathbb{Q}_{\geq 0}$, find a maximum cut $(X',Y')$ of
$F_k(G)$.
\end{problem}

\begin{problem}[\textsc{Token Max-Cut}, decision version]\label{prob:tmc_dec}
Given $F_k(G)$ and $w$ as above, and a rational $t$,
decide whether $m(F_k(G))\geq t$.
\end{problem}

We allow the weights to be arbitrary, as this is the form in which the phase
reconstruction problem of the introduction presents them.
Problem~\ref{prob:tmc} asks for a cut and not merely for its weight, as it is the
cut itself that prescribes the signs of the ground state in the phase
reconstruction problem; Problem~\ref{prob:tmc_dec} is its
decision version.
Problem~\ref{prob:tmc_dec} is in NP because a cut of $F_k(G)$ has size polynomial in $N$.

%

\begin{theorem}\label{thm:tmc_npc}
For every constant $k$, Problem~\ref{prob:tmc_dec} is  NP-hard, when
restricted to unweighted instances.
\end{theorem}

\begin{proof}

We use the following result on the hardness of approximating
\textsc{Max-Cut}~\cite{Hastad01, TSSW00, crescenzi2001weighted}.
For every constant $16/17<\rho<1$, there is a polynomial-time algorithm that, given a
$3$-\textsc{Sat} formula $\varphi$, outputs an unweighted graph $G_\varphi$ and a positive
integer $t$ such that:
\begin{itemize}
\item if $\varphi$ is satisfiable, then $m(G_\varphi)\geq t$; and
 \item if $\varphi$ is not satisfiable, then $m(G_\varphi)<\rho\, t$.
\end{itemize}
Fix one such $\rho$.

By $(a)$ of Theorem~\ref{thm:upper_bound} and since $k$ is constant, there exists a natural number $n_0$
such that if $G$ has $n \ge n_0$ vertices, then
\[m(F_k(G)) \le \frac{1}{\rho}\binom{n-2}{k-1}m(G).\]
Since adding isolated vertices to $G_\varphi$ does not change the value
of $m(G_\varphi)$, we may assume that $G_\varphi$ has $n \ge n_0$ vertices. Let
\[t':=\binom{n-2}{k-1}t.\]
Since $k$ is constant, $F_k(G_\varphi)$ has $O(n^k)$ vertices; thus, $F_k(G_\varphi)$ and $t'$
can be computed from $\varphi$ in polynomial time. We have that
\begin{itemize}
\item if $\varphi$ is satisfiable, then by (\ref{eq:lower}),
\[m(F_k(G_\varphi)) \ge \binom{n-2}{k-1}m(G_\varphi) \ge \binom{n-2}{k-1} t=t'; \]
\item if $\varphi$ is not satisfiable, then
\[m(F_k(G_\varphi)) \le \frac{1}{\rho}\binom{n-2}{k-1}m(G_\varphi) < \frac{1}{\rho}\binom{n-2}{k-1}\rho\, t=t'. \]
\end{itemize}
Therefore, $\varphi$ is satisfiable if and only if $m(F_k(G_\varphi)) \ge t'$; and the result follows.
\end{proof}

We now give evidence that  for every $\delta >0$
and $k \ge n^\delta$, Problem~\ref{prob:tmc_dec} is \emph{not} NP-hard,  provided the weights are induced by a weighting of $G$
taking finitely many values. The argument  says nothing when the weights of $F_k(G)$ are arbitrary.

A Turing machine \emph{with advice} receives, besides its input $x$, a second
string $\alpha_{|x|}$ that depends only on the length of $x$; apart from this
it is an ordinary deterministic machine. The advice is not required to be
computable, and the same string must serve every input of that length. Let
$\mathrm{P}/f(n)$ be the set of languages decided in polynomial time with
advice of length $f(n)$; the most studied of these classes is
\[
\mathrm{P}/\mathrm{poly}:=\bigcup_{c\geq 1}\mathrm{P}/n^{c},
\]
the languages decided in polynomial time with advice of polynomial length.
Clearly
$\mathrm{P}\subseteq\mathrm{P}/\mathrm{poly}$.
It is widely believed that $\mathrm{NP}\not\subseteq\mathrm{P}/\mathrm{poly}$.
The evidence comes from the following theorem of Karp and Lipton~\cite{KarpLipton1980}. If
$\mathrm{NP} \subseteq \mathrm{P}/\mathrm{poly}$, then the polynomial hierarchy
would collapse to its second level. We refer the reader to Chapters~5 and~6
of~\cite{AroraBarak2009} for the polynomial hierarchy, for advice classes and
circuits, and for the theorem of Karp and Lipton.

We say that $f(n)$ is \emph{quasipolynomial} if
\[
f(n)=2^{O\left((\log n)^{c}\right)}
\]
for some constant $c>0$.
Let $W$ be a finite set of nonnegative weights, fixed once and for all. We call
an instance $\bigl(F_k(G),w,t\bigr)$ of Problem~\ref{prob:tmc_dec}
\emph{$W$-weighted} if $w$ is induced by a weighting of $G$ with values in $W$;
the unweighted case is that of $W=\{1\}$.
For $k\geq n^{\delta}$ we can use advice of quasipolynomial length to solve
the $W$-weighted case of
Problem~\ref{prob:tmc_dec} in quasipolynomial time. The main idea is that, for a
fixed input size, there is only a quasipolynomial number of weighted parent
graphs that could have given rise to $F_k(G)$.

\begin{prop}\label{prop:advice}
Let $\delta>0$ be fixed and let $k=k(n)\geq n^{\delta}$. Then for this value of
$k$, the $W$-weighted case of Problem~\ref{prob:tmc_dec} can be solved in
quasipolynomial time with advice of quasipolynomial length.
\end{prop}

\begin{proof}
Let $\bigl(F_k(G),w,t\bigr)$ be a $W$-weighted instance of  Problem~\ref{prob:tmc_dec}, and let $s$ be its size. As $F_k(G)$
has $\binom{n}{k}$ vertices and $k\leq n/2$,
\[
s\;\geq\;\binom{n}{k}\;\geq\;\left(\frac{n}{k}\right)^{k}\;\geq\;2^{k}\;\geq\;2^{n^{\delta}},
\]
so that $n\leq(\log s)^{1/\delta}$.

As advice for the inputs of size $s$ we take a table of the values
$m\bigl(F_k(G)\bigr)$, over all $W$-weighted labeled graphs $G$ on at most
$(\log s)^{1/\delta}$ vertices. There are at most
$(|W|+1)^{\binom{n}{2}}=2^{O\left((\log s)^{2/\delta}\right)}$ such graphs, and
each value is a sum of at most $s$ elements of $W$, and can thus be written with
$O(\log s)$ bits. So the table has length quasipolynomial in $s$; and,
as is required of advice, it depends only on $s$.

Given the instance, we recover $G$ and its weighting from the labels of the
vertices of $F_k(G)$, which takes polynomial time since $ab\in E(G)$ precisely
when $S\cup\{a\}$ and $S\cup\{b\}$ are adjacent in $F_k(G)$ for some $S$, and
the weight of $ab$ is that of any such edge. We then read
$m\bigl(F_k(G)\bigr)$ off the table and compare it with $t$. Locating the entry
of $G$ takes time proportional to the length of the table, so the whole computation is quasipolynomial. 
\end{proof}

Note that in Proposition~\ref{prop:advice}, time is quasipolynomial in terms of $N$. For example for $k=n/2$, we have that $N=\Theta(2^n/\sqrt{n})$, which is already exponential on $n$.

Suppose now that the $W$-weighted case of Problem~\ref{prob:tmc_dec}, restricted
to some $k\geq n^{\delta}$, were NP-hard. Composing a reduction from \textsc{Sat} with the
algorithm of Proposition~\ref{prop:advice} would then decide an NP-complete set
in quasipolynomial time with advice of quasipolynomial length. This is the
quasipolynomial counterpart of the hypothesis
$\mathrm{NP}\subseteq\mathrm{P}/\mathrm{poly}$ discussed above, and it is
believed to be false for an analogous reason: by a theorem of Buhrman and
Homer~\cite[Theorem~1]{BuhrmanHomer1992}, it would collapse the exponential
hierarchy to
its second level, $\mathrm{NEXP}^{\mathrm{NP}}$, just as
$\mathrm{NP}\subseteq\mathrm{P}/\mathrm{poly}$ collapses the polynomial
hierarchy to its second level.
\section*{Declaration of AI use}

Large language models were used to draft portions of the text, to search the literature, to verify small cases numerically, and for copy-editing; a number of the proof ideas were obtained in discussions with them. Everything has been checked by the authors, who are solely responsible for the correctness of the final text.

\bibliographystyle{alphaurl}

\bibliography{maxcut_new}

\end{document}